\documentclass[12pt]{article}
\usepackage{graphicx}
\scrollmode
\usepackage{amsfonts}
\usepackage{amsmath}
\usepackage{amssymb}
\begin{document}
\title{\bf The Brownian Web}
\author{L.~R.~G.~Fontes {\thanks {Instituto de 
Matem\'{a}tica e Estatistica, Universidade de S\~{a}o Paulo,
05311-970 S\~{a}o Paulo SP, Brasil}}
\and M.~Isopi {\thanks {Dipartimento di Interuniversitario Matematica,
Universit\`{a} di Bari, 70125 Bari, Italia}}\and C.~M.~Newman 
{\thanks {Courant Institute of Mathematical Sciences, New York University,
New York, NY 10012}}\and K.~Ravishankar {\thanks {Department of Mathematics,
State University of New York, College at New Paltz, New Paltz, NY 12561}}}
\maketitle

\newtheorem{defin}{Definition}[section]
\newtheorem{Prop}{Proposition}
\newtheorem{teo}{Theorem}[section]
\newtheorem{ml}{Main Lemma}
\newtheorem{con}{Conjecture}
\newtheorem{cond}{Condition}
\newtheorem{prop}[teo]{Proposition}
\newtheorem{lem}{Lemma}[section]
\newtheorem{rmk}[teo]{Remark}
\newtheorem{cor}{Corollary}[section]
\renewcommand{\theequation}{\thesection .\arabic{equation}}

\newcommand{\beq}{\begin{equation}}
\newcommand{\eeq}{\end{equation}}
\newcommand{\beqn}{\begin{eqnarray}}
\newcommand{\beqnn}{\begin{eqnarray*}}
\newcommand{\eeqn}{\end{eqnarray}}
\newcommand{\eeqnn}{\end{eqnarray*}}
\newcommand{\bprop}{\begin{prop}}
\newcommand{\eprop}{\end{prop}}
\newcommand{\bteo}{\begin{teo}}
\newcommand{\bcor}{\begin{cor}}
\newcommand{\ecor}{\end{cor}}
\newcommand{\bcon}{\begin{con}}
\newcommand{\econ}{\end{con}}
\newcommand{\bcond}{\begin{cond}}
\newcommand{\econd}{\end{cond}}
\newcommand{\eteo}{\end{teo}}
\newcommand{\brm}{\begin{rmk}}
\newcommand{\erm}{\end{rmk}}
\newcommand{\blem}{\begin{lem}}
\newcommand{\elem}{\end{lem}}
\newcommand{\ben}{\begin{enumerate}}
\newcommand{\een}{\end{enumerate}}
\newcommand{\bei}{\begin{itemize}}
\newcommand{\eei}{\end{itemize}}
\newcommand{\bdf}{\begin{defin}}
\newcommand{\edf}{\end{defin}}

\newcommand{\nn}{\nonumber}
\renewcommand{\=}{&=&}
\renewcommand{\>}{&>&}
\renewcommand{\le}{\leq}
\newcommand{\+}{&+&}
\newcommand{\fr}{\frac}

\renewcommand{\r}{{\mathbb R}}
\newcommand{\br}{\bar{\mathbb R}}
\newcommand{\Z}{{\mathbb Z}}
\newcommand{\z}{{\mathbb Z}}
\newcommand{\zd}{\z^d}
\newcommand{\zz}{{\mathbb Z}}
\newcommand{\R}{{\mathbb R}}
\newcommand{\tw}{\tilde{\cal W}}
\newcommand{\E}{{\mathbb E}}
\newcommand{\C}{{\mathbb C}}
\renewcommand{\P}{{\mathbb P}}
\newcommand{\N}{{\mathbb N}}
\newcommand{\var}{{\mathbb V}}
\renewcommand{\S}{{\cal S}}
\newcommand{\T}{{\cal T}}
\newcommand{\cm}{{\cal M}}
\newcommand{\cp}{{\cal P}}
\newcommand{\h}{{\cal H}}
\newcommand{\f}{{\cal F}}
\newcommand{\cd}{{\cal D}}
\newcommand{\xt}{X_t}
\renewcommand{\ge}{g^{(\epsilon)}}
\newcommand{\xe}{y^{(\epsilon)}}
\newcommand{\ye}{y^{(\epsilon)}}
\newcommand{\bx}{{\bar y}}
\newcommand{\by}{{\bar y}}
\newcommand{\bxe}{{\bar y}^{(\epsilon)}}
\newcommand{\bye}{{\bar y}^{(\epsilon)}}
\newcommand{\bwe}{{\bar w}^{(\epsilon)}}
\newcommand{\bxz}{{\bar y}}
\newcommand{\bwz}{{\bar w}}
\newcommand{\we}{w^{(\epsilon)}}
\newcommand{\Xe}{Y^{(\epsilon)}}
\newcommand{\Ze}{Z^{(\epsilon)}}
\newcommand{\Ye}{Y^{(\epsilon)}}
\newcommand{\ydo}{Y^{(\d)}_{y_0(\d),s_0(\d)}}
\newcommand{\yo}{Y_{y_0,s_0}}
\newcommand{\tye}{{\tilde Y}^{(\epsilon)}}
\newcommand{\hy}{{\hat Y}}
\newcommand{\ve}{V^{(\epsilon)}}
\newcommand{\Ne}{N^{(\epsilon)}}
\newcommand{\ce}{c^{(\epsilon)}}
\newcommand{\cle}{c^{(\l\epsilon)}}
\newcommand{\xet}{Y^{(\epsilon)}_t}
\newcommand{\hxt}{\hat X_t}
\newcommand{\btn}{\bar\tau_n}
\newcommand{\ct}{{\cal T}}
\newcommand{\rn}{{\cal R}_n}
\newcommand{\nt}{{N}_t}
\newcommand{\lnk}{{\cal L}_{n,k}}
\newcommand{\cl}{{\cal L}}
\newcommand{\bw}{\bar{\cal W}}
\newcommand{\tc}{\tilde{\cal C}_b}
\newcommand{\hxtt}{\hat X_{\ct}}
\newcommand{\txnt}{\tilde X_{\nt}}
\newcommand{\xs}{X_s}
\newcommand{\xn}{\tilde X_n}
\newcommand{\tx}{\tilde X}
\newcommand{\hx}{\hat X}
\newcommand{\txi}{\tilde X_i}
\newcommand{\txij}{\tilde X_{i_j}}
\newcommand{\taxi}{\tau_{\txi}}
\newcommand{\txn}{\tilde X_N}
\newcommand{\xk}{X_K}
\newcommand{\ts}{\tilde S}
\newcommand{\tl}{\tilde\l}
\newcommand{\tg}{\tilde g}
\newcommand{\im}{I^-}
\newcommand{\ip}{I^+}
\newcommand{\hal}{H_\a}
\newcommand{\ba}{B_\a}

\renewcommand{\a}{\alpha}
\renewcommand{\b}{\beta}
\newcommand{\g}{\gamma}
\newcommand{\G}{\Gamma}
\renewcommand{\L}{\Lambda}
\renewcommand{\d}{\delta}
\newcommand{\D}{\Delta}
\newcommand{\e}{\epsilon}
\newcommand{\fes}{\phi^{(\epsilon)}_s}
\newcommand{\fet}{\phi^{(\epsilon)}_t}
\newcommand{\fe}{\phi^{(\epsilon)}}
\newcommand{\pset}{\psi^{(\epsilon)}_t}
\newcommand{\pse}{\psi^{(\epsilon)}}
\renewcommand{\l}{\lambda}
\newcommand{\me}{\mu^{(\epsilon)}}
\newcommand{\re}{\rho^{(\epsilon)}}
\newcommand{\tre}{\tilde{\rho}^{(\epsilon)}}
\newcommand{\nue}{\nu^{(\epsilon)}}
\newcommand{\mbe}{{\bar\mu}^{(\epsilon)}}
\newcommand{\rbe}{{\bar\rho}^{(\epsilon)}}
\newcommand{\mb}{{\bar\mu}}
\newcommand{\rb}{{\bar\rho}}
\newcommand{\mbz}{{\bar\mu}}
\newcommand{\s}{\sigma}
\renewcommand{\o}{\Pi}
\newcommand{\om}{\omega}
\newcommand{\tio}{\tilde\o}
\renewcommand{\sl}{\sigma'}
\newcommand{\si}{\s(i)}
\newcommand{\sit}{\s_t(i)}
\newcommand{\ei}{\eta(i)}
\newcommand{\eit}{\eta_t(i)}
\newcommand{\eot}{\eta_t(0)}
\newcommand{\sil}{\s'_i}
\newcommand{\sj}{\s(j)}
\newcommand{\st}{\s_t}
\newcommand{\so}{\s_0}
\newcommand{\xii}{\xi_i}
\newcommand{\xij}{\xi_j}
\newcommand{\xio}{\xi_0}
\newcommand{\ti}{\tau_i}
\newcommand{\te}{\tau^{(\epsilon)}}
\newcommand{\bt}{\bar\tau}
\newcommand{\tti}{\tilde\tau_i}
\newcommand{\tto}{\tilde\tau_0}
\newcommand{\tei}{T_i}
\newcommand{\ttei}{\tilde T_i}
\newcommand{\tes}{T_S}
\newcommand{\tao}{\tau_0}

\renewcommand{\t}{\tilde t}

\newcommand{\da}{\downarrow}
\newcommand{\ua}{\uparrow}
\newcommand{\ar}{\rightarrow}
\newcommand{\lar}{\leftrightarrow}
\newcommand{\va}{\stackrel{v}{\rightarrow}}
\newcommand{\ppa}{\stackrel{pp}{\rightarrow}}
\newcommand{\dw}{\stackrel{w}{\Rightarrow}}
\newcommand{\Va}{\stackrel{v}{\Rightarrow}}
\newcommand{\Ppa}{\stackrel{pp}{\Rightarrow}}
\newcommand{\la}{\langle}
\newcommand{\ra}{\rangle}
\newcommand{\ep}{\vspace{.5cm}}
\newcommand\sqr{\vcenter{
           \hrule height.1mm
           \hbox{\vrule width.1mm height2.2mm\kern2.18mm\vrule width.1mm}
           \hrule height.1mm}}                  


\begin{abstract}
Arratia, and later T\'oth and Werner, constructed random processes
that formally correspond to coalescing one-dimensional Brownian motions starting
from every space-time point. We 
extend their work
by constructing and characterizing
what we call the {\em Brownian Web} as a random variable taking values in
an appropriate (metric) space whose points are (compact) sets of paths.
This leads to general convergence criteria and, in particular,
to convergence in distribution of coalescing random walks 
in the scaling limit to the Brownian 
Web.
\end{abstract}

\vspace{.75cm}


\noindent{\Large{\bf Introduction}}

\vspace{.5cm}

Construct random paths in the plane, as follows.
Take the square lattice consisting of all points $(\sqrt{2}m,\sqrt{2}n)$
with $m,n$ integers and rotate it by $45$ degrees resulting in
all points $(i,j)$ with $i,j$ integers and $i+j$ even.
Imagine a walker at spatial location $i$ at time $j$ deciding
to move right or left at unit speed between times $j$ and $j+1$
if the outcome of a fair coin toss is heads ($ \Delta_{i,j} = +1 $) 
or tails ($ \Delta_{i,j} = -1 $), with the coin tosses independent
for different space-time points $(i,j)$. Figure 1 depicts a
simulation of the resulting paths.

The path of a walker starting from $y_0$ at time $s_0$ is the graph
of a simple 
symmetric one-dimensional random walk, $\yo (t)$. At integer times, $\yo (t)$
is the solution of the simple stochastic difference equation,
\begin{equation}
\label{diffeq}
Y(j+1) - Y(j) =  \Delta_{Y(j),j} , \quad Y(s_0) = y_0.
\end{equation}
Note that the paths of distinct walkers starting from different
$(y_0,s_0)$'s are automatically {\em coalescing} --- i.e., they are 
independent of each other
until they coalesce (i.e., become identical) upon meeting
at some space-time point.
  
After rescaling to spatial steps of size $\d$ and time steps
of size $\d^2$, a single rescaled random walk (say, starting from
$0$ at time $0$) $Y_{0,0}^{(\d)}(t) = \d Y_{0,0}(\d^{-2}t)$
converges as $\d \to 0$ to a standard Brownian motion $B(t)$. More precisely,
by the Donsker invariance principle~\cite{kn:D},
the distribution of $Y_{0,0}^{(\d)}$ on the space of continuous paths
converges weakly as $\d\to0$ to standard Wiener measure.

The invariance principle is also valid for continuous time random walks,
where the move from $i$ to $i\pm1$ takes an exponentially distributed
time (see the discussion following Remark~\ref{rmk:etaconv}
below for more details). 
In continuous time, coalescing random walks are at the
heart of Harris's graphical representation
of the (one-dimensional) voter model~\cite{kn:H} and their scaling limits arise
naturally in the physical context of (one-dimensional) aging~\cite{kn:FINS}.
Like for a single random walk, finitely many rescaled coalescing walks
in discrete or continuous time (with rescaled space-time starting points) 
converge in distribution to  finitely many coalescing
Brownian motions. In this paper, we present results
concerning the convergence
in distribution of the collection of the rescaled coalescing walks from
{\em all} the starting points; 
detailed proofs will be published elsewhere~\cite{kn:FINR}.

\begin{figure}[!ht]
\begin{center}
\includegraphics[width=10cm]{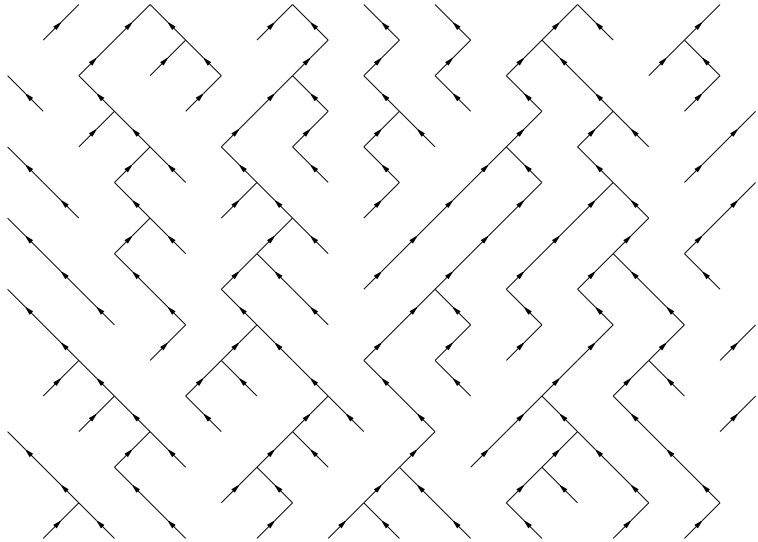}
\caption{Coalescing random walks in discrete time; 
the horizontal coordinate is space and the 
vertical one is time.}\label{figpath}
\label{path} 
\end{center}
\end{figure}

Our results come in two parts:
\begin{itemize}
\item[(1)] characterization (and construction) of the limiting object,
which we call the standard {\em Brownian Web (BW)}, and
\item[(2)] general convergence criteria, which are then applied to
coalescing random walks.
\end{itemize}
A key ingredient of the characterization and construction
(see Theorem~\ref{thm:char}) 
is the choice of a space for the Brownian web; this is the
BW analogue of the space of continuous paths for Brownian motion.
The convergence criteria and application (see
Theorems~\ref{thm:conv} and~\ref{thm:crwtobw} below)
are the BW analogues of Donsker's invariance principle.
Like Brownian motion itself, 
we expect that the Brownian web
and its variants (see, e.g., Remark~\ref{rmk:dual}) will be quite 
ubiquitous as scaling limits, well beyond the context of 
coalescing random walks and our sufficient conditions for convergence.

Much of the construction of the Brownian web
was already done in the groundbreaking work of
Arratia~\cite{kn:A1,kn:A2} and then in work of T\' oth and
Werner~\cite{kn:TW} (see also~\cite{kn:STW}). They 
all recognized that in the
limit $\d\to0$ there would be (nondeterministic) space-time points
$(x,t)$ starting from which 
there are multiple limit paths and they provided
various conventions (e.g., semicontinuity in $x$) to avoid such
multiplicity. Our main contribution vis-a-vis
construction is to accept the intrinsic
nonuniqueness by choosing an appropriate metric space in which the BW takes
its values. Roughly speaking, instead of using some convention to obtain
a process that is a {\em single-valued} mapping from 
each space-time starting
point to a single path from that starting point, we allow 
{\em multi-valued} mappings; more accurately, our BW value is the collection
of {\em all} paths from all starting points. This choice of space is 
very much in the spirit of earlier
work~\cite{kn:A,kn:AB,kn:ABNW} on spatial scaling limits of
critical percolation models and
spanning trees, but modified for our particular space-time setting.


\section{Brownian Web: Characterization}
\setcounter{equation}{0}
\label{sec:cons}

We begin by defining three metric spaces: $(\br^2,\rho)$, $(\o,d)$ and 
$(\h,d_\h)$. 
The elements of the three spaces are respectively: points
in space-time, paths with specified starting points in space-time and
collections of paths with specified starting points. 
The BW will be an $(\h,\f_\h)$-valued random variable,
where $\f_\h$ is the Borel $\s$-field  associated to the metric $d_\h$.

$(\br^2,\rho)$ is the completion (or compactification) of $\R^2$ under the
metric $\rho$, where
\begin{equation}
\label{rho}
\rho((x_1,t_1),(x_2,t_2))=
\left|\frac{\tanh(x_1)}{1+|t_1|}-\frac{\tanh(x_2)}{1+|t_2|}\right|
\vee|\tanh(t_1)-\tanh(t_2)|.
\end{equation}
$\br^2$ may be thought as the set of $(x,t)$ in 
$[-\infty,\infty]\times[-\infty,\infty]$ with all points of the form
$(x,-\infty)$ identified (and similarly for $(x,\infty)$). More
precisely, it is the image of $[-\infty,\infty]\times[-\infty,\infty]$
under the mapping
\begin{equation}
\label{compactify}
(x,t)\leadsto(\Phi(x,t),\Psi(t))
\equiv\left(\frac{\tanh(x)}{1+|t|},\tanh(t)\right).
\end{equation}

For $t_0\in[-\infty,\infty]$, let $C[t_0]$ denote the set of functions
$f$ from $[t_0,\infty]$ to $[-\infty,\infty]$ such that $\Phi(f(t),t)$
is continuous. Then define 
\begin{equation}
\label{omega}
\o=\bigcup_{t_0\in[-\infty,\infty]}C[t_0]\times\{t_0\},
\end{equation}
where $(f,t_0)\in\o$ then represents a path in $\br^2$ starting at $(f(t_0),t_0)$.
For$(f,t_0)\in\o$, we denote by $\hat f$ the function that extends $f$ to all
$[-\infty,\infty]$ by setting it equal to $f(t_0)$ for $t<t_0$. Then we take
\begin{equation}
\label{d}
d((f_1,t_1),(f_2,t_2))=
(\sup_t|\Phi(\hat{f_1}(t),t)-\Phi(\hat{f_2}(t),t)|)
\vee|\Psi(t_1)-\Psi(t_2)|.
\end{equation}
$(\o,d)$ is a complete separable metric space. 

Let now $\h$ denote the set of compact
subsets of $(\o,d)$, with $d_\h$ the induced Hausdorff metric, i.e.,
\begin{equation}
\label{dh}
d_\h(K_1,K_2)=\sup_{g_1\in K_1}\inf_{g_2\in K_2}d(g_1,g_2)\vee
              \sup_{g_2\in K_2}\inf_{g_1\in K_1}d(g_1,g_2).
\end{equation}
$(\h,d_\h)$ is also a complete separable metric space.

Before stating our characterization theorem for the Brownian web, we need 
some definitions. For an $(\h,\f_\h)$-valued random variable $\bar W$ (or its
distribution $\mu$), we define the {\em finite-dimensional distributions}
of $\bar W$ as the induced probability measures $\mu_{(x_1,t_1;\ldots;x_n,t_n)}$
on the subsets of paths starting from any 
finite deterministic set of points
$(x_1,t_1),\ldots,(x_n,t_n)$ in $\R^2$.
There are several ways in which the Brownian web can be
characterized; they differ from each other primarily
in the type of extra condition required
beyond the
finite-dimensional distributions. 
The characterization of the next theorem, or more precisely a variant
discussed later in Remark~\ref{rmk:cons}, is the one most directly
suited to the convergence results of Section~\ref{sec:conv};
an alternative characterization in which the extra condition
is a type of Doob separability property (see, e.g., Chap.~3 of \cite{kn:V})
is discussed 
in Remark~\ref{rmk:separable}. For
the next theorem, we also define, for $t\geq0$ and $a\leq b$, the 
($\{0,1,\ldots,\infty\}$-valued) random variable
$\eta(t_0,t;a,b)$ as the number of {\em distinct} points in
$\R\times\{t_0+t\}$ that are touched by paths in $\bar W$ which also touch 
some point in $[a,b]\times\{t_0\}$.

\bteo
\label{thm:char}
There is an $(\h,\f_\h)$-valued random variable $\bar W$ whose distribution
$\mu$ is uniquely determined by the following two properties:
\begin{itemize}
\item[(i)] its finite-dimensional distributions are those of coalescing
Brownian motions (with unit diffusion constant), and
\item[(ii)] for $-\infty<t_0<\infty$, $0<t<\infty$, $-\infty<a\leq b<\infty$,
\begin{equation}
\label{eta}
E_\mu(\eta(t_0,t;a,b))= 1 + \frac{b-a}{\sqrt{\pi t}}.
\end{equation}
\end{itemize}

\eteo

\brm
\label{rmk:separable}
Implicit in condition~(i) of the theorem is that starting from
any deterministic point, there is almost
surely only a single path in $\bar W$.
Condition (ii) can be replaced by the separability property
that there is a deterministic dense countable set $\cd$ of
space-time starting points, such that
almost surely, $\bar W$ is the closure in $(\o,d)$ of the set of
paths starting from the points of $\cd$. It should be noted
that $(\h,\f_\h)$-valued random variables,
satisfying condition (i) but not condition (ii) or its separability
alternative, can occur naturally.
Such a process (closely related to the ``Double Brownian Web'' of 
Remark~\ref{rmk:dual} below), where the counting variable $\eta$
is infinite with strictly positive probability, will be studied
elsewhere and shown to arise as the scaling limit of stochastic flows,
extending earlier work of Piterbarg \cite{kn:P}.
\erm

\noindent{\bf Sketch of Proof of Theorem~\ref{thm:char}.} 
The construction of the Brownian web
(i.e., the existence of such a $\bar W$) begins as
in ~\cite{kn:A2,kn:TW} with the construction of a set $W$ of
coalescing Brownian paths starting from a deterministic dense countable set
$\cd$ of space-time starting points. This skeleton 
$W = \{\tilde{W}_1,\tilde{W}_2,\dots,\}$
is a random subset of $\o$ that is constructed by 
deterministically ordering the points
of $\cd$ as $(x_1,t_1), (x_2,t_2), \dots$, then defining
$W_j = (x_j+B_j(t-t_j),t_j)\in \o$ where the $B_j$'s are independent
standard Brownian motions, and finally using the ordering to inductively define
$\tilde{W}_j \in \o$ by following $W_j$ until it meets some $\tilde{W}_k$
with $k<j$ after which point it follows $\tilde{W}_k$. 

The next several steps of the construction are to show that the 
closure $\bar W$
in $(\o,d)$ of this BW-skeleton is compact, that the distribution of
$\bar W$  does not depend on the choice of $\cd$ or its ordering, 
and that $\bar W$
satisfies (i) and (ii) of Theorem~\ref{thm:char} above. The compactness
can be proved in a number of ways; one of these is to verify a 
condition, as in (\ref{T1}) below, but with $\mu_\d$ replaced by
the distribution of $\{\tilde{W}_1,\dots,\tilde{W}_m\}$ and the sup
over $\d$ replaced by a sup over $m$ (and then argue as 
at the beginning of the proof of Theorem~\ref{thm:conv} below, eventually
invoking the Arzel\`a-Ascoli theorem). To verify the said condition, 
one argues
as in the last two paragraphs of the proof of 
Theorem~\ref{thm:crwtobw} below. The 
argument actually involves only a single bound
like~(\ref{gbound1}), which is obtained in the same way as 
in the proof of 
Theorem~\ref{thm:crwtobw}.
We remark that by considering the quantity $\tg(t,u)$, as in (\ref{T1}),
but with $u = t^{\xi}$, one can show not just compactness
of $\bar W$, but also H\"{o}lder continuity with any exponent $\xi < 1/2$
for {\em all} the paths of $\bar W$. 

The lack of dependence of the distribution
on the choice of $\cd$ or its ordering
follows fairly directly after verifying property (i)
for $\bar W$. Property (i) itself follows by a trapping argument about 
a deterministic point $(\bar x,\bar t)$ (and similarly for finitely 
many points)
and any sequence $(\bar{x}_i,\bar{t}_i)=(x_{j(i)},t_{j(i)})$ 
from $\cd$ converging to $(\bar{x},\bar{t})$ as 
$i\to\infty$: 
that (even if $j(i)$ is nondeterministic and regardless
of whether $t_{j(i)}>\bar{t}$ or $t_{j(i)} \le \bar{t}$) for large $i$,
$\tilde{W}_{j(i)}$ is 
(with probability very close to one) trapped between 
$\tilde{W}_k$ and $\tilde{W}_{k'}$ for deterministic $k,k'$ with $x_k<\bar{x}<x_k'$,
$t_k<\bar t$, $t_{k'}<\bar t$, $x_k$ and $x_k'$ close to $\bar x$ 
and $t_k,t_{k'}$ even closer to $\bar t$
so that $\tilde{W}_{j(i)}$ (with probability very close to one) quickly coalesces with
both $\tilde{W}_k$ and $\tilde{W}_{k'}$; thus $\tilde{W}_{j(i)}$ converges
almost surely as $i \to \infty$ to 
a path (independent of the specific sequence $(\bar{x}_i,\bar{t}_i)$)
that is distributed as a Brownian motion starting from $(\bar{x},\bar{t})$.

Verifying property (ii) is somewhat indirect. First, one shows that the 
random variable $\eta$ for our constructed $\bar W$ is almost surely finite
with a finite mean which, by the translation invariance in space and time
that results from the lack of dependence on $\cd$, must be of the form
$\Lambda(b-a,t)$. Second, the specific evaluation of $\Lambda$ as given on the righthand
side of (\ref{eta}) is carried out.  As the explicit expression for $\Lambda$ will not
actually be used in our convergence results of the next section (see Remark~\ref{rmk:cons}),
we can and will use those convergence 
results in the evaluation of $\Lambda$. 

The first part of the verification of (ii) is a consequence of an inequality,
\begin{equation}
\label{etabound}
P(\eta(t_0,t;a,b) \geq k) \le [P(\eta(t_0,t;a,b) \geq 2)]^{k-1} = 
[\Theta(b-a,t)]^{k-1},
\end{equation}
where $\Theta(b-a,t)$ is the 
probability that two independent Brownian motions
starting at a distance $b-a$ apart at 
time zero will have met by time $t$ (which
itself can be expressed in terms of a single Brownian motion). 
The inequality in (\ref{etabound}) is first derived for finite subsets
$\{\tilde{W}_1,\tilde{W}_2,\dots,\tilde{W}_m\}$ of the 
skeleton, and thus for the whole skeleton. For
the whole skeleton and its closure $\bar W$, the equality 
in (\ref{etabound}) is seen to be valid by
choosing the countable set $\cd$ so that its first two points 
are $(a,t_0)$ and $(b,t_0)$.
Then the inequality is extended to $\bar W$ from the 
skeleton by a limit/approximation
argument which uses that $\{K\in \h\,:\, \tilde{\eta}(K) \geq k \}$
is open in $(\h,d_\h)$, where $\tilde{\eta}$ is the modification 
of $\eta$ that counts points in $\R\times\{t_0+t\}$ touched by paths
in $K$ which touch 
$(a-\tilde{\varepsilon}_1,b+\tilde{\varepsilon}_1) 
\times \{t_0+\tilde{\varepsilon}_2\}$
and start earlier than $t_0+\tilde{\varepsilon}_2$. 

The second part of the 
verification of (ii), in which $\Lambda$ is explicitly evaluated,
is a consequence of all the following: 
a result of Bramson and Griffeath~\cite{kn:BG} on the large-time
asymptotics of mean
interparticle distance in coalescing random walks, the conversion
of that result by standard arguments to asymptotics for the 
mean of the rescaled random walk version of the counting
variable $\eta$, convergence of the distribution of $\eta$
in the scaling limit (see Remark~\ref{rmk:etaconv}), and finally
the analogue of (\ref{etabound}) 
for coalescing walks (see (\ref{walkbound}))
which implies uniform integrability of $\eta$
as $\d \to 0$ and hence convergence of the mean of $\eta$.

It remains to show that conditions 
(i) and (ii) for a measure $\mu'$ on $(\h,\f_\h)$
together imply that $\mu'$ equals the 
distribution $\mu$ of the constructed Brownian web $\bar W$. Let us denote by 
$X'$ the $(\h,\f_\h)$-valued random variable distributed by $\mu'$ and 
by $\eta'$ the counting random variable appearing in condition (ii) for $\mu'$.
Choose some deterministic dense countable subset $\cd$ and consider the countable
collection $W^*$ of paths of $X'$ starting from $\cd$. By condition (i), $W^*$
is equidistributed with our constructed 
Brownian skeleton $W$ (based on the same $\cd$)
and hence the closure $\bar{W^*}$ of $W^*$ in $(\o,d)$ is a subset of $X'$
that is equidistributed with our constructed Brownian
web $\bar W$. To complete the proof, we will use condition (ii) to show
that $X' \setminus \bar{W^*}$ is almost surely empty
by using the fact that 
the counting variable $\eta^*$ for $\bar{W^*}$ already satisifies
condition (ii) since $\bar{W^*}$ is distributed as a Brownian web. If
$X' \setminus \bar{W^*}$ were nonempty (with strictly positive probability),
then there would have to be some rational $t_0, t, a, b$ 
for which $ \eta' > \eta^*$. But then
\begin{equation}
P(\eta'(t_0,t;a,b) > \eta^*(t_0,t;a,b)) \,>\, 0 
\end{equation}
for some rational $t_0, t, a, b$ and 
so condition (ii) for $\eta'$ would 
{\em not} be valid for that $t_0, t, a, b$.  

\brm
\label{rmk:cons}
The proof of Theorem~\ref{thm:char} makes clear that
the idea behind (i) and (ii) together implying uniqueness of the distribution
is that (i) implies sufficiently many paths and (ii) implies no extraneous
ones. 
Thus condition (i) can be weakened to the existence of a subset of paths
distributed as the coalescing Brownian motions of the skeleton $W$ 
(for any deterministic dense countable $\cd$) and
condition (ii) can also be modified, e.g., by replacing the equality
in~(\ref{eta}) by an inequality ($\leq$) 
and by replacing an (in)equality for the mean
by one for the distribution. Similarly, in applying our characterization
results to obtain convergence criteria as we do in  Theorem~\ref{thm:conv},
an explicit expression for the mean as given in the righthand 
side of (\ref{eta}) or an explicit expression for the distribution
is not needed; i.e., to verify that an $X'$ is equidistributed
with our explicitly constructed Brownian web $\bar W$,
condition (ii) for the $\eta'$ of $X'$ can 
be replaced by the condition
that the distribution of $\eta'$ equal (or only is 
stochastically dominated by)
the distribution of the $\eta$ of $\bar W$. 
\erm

\brm
\label{rmk:dual}
In the graphical representation of Harris
for the one-dimensional voter model~\cite{kn:H},
coalescing random walks forward in time and coalescing dual random walks
backward in time (with forward and backward walks not crossing each other)
are constructed simultaneously (see, e.g., the discussion in~\cite{kn:FINS}).
The simultaneous construction of forward and (dual) backward Brownian motions
was emphasized in~\cite{kn:TW,kn:STW} and their approach and results can be 
applied to extend both our characterization and convergence results to the
{\em Double Brownian Web} (DBW) which includes simultaneously the forward BW
and its dual backward BW. We note that in the DBW, the $\eta$ of (\ref{eta})  
equals $1+ \eta^{\mbox{\rm\scriptsize dual}}$, where 
$\eta^{\mbox{\rm\scriptsize dual}}$ is the number of distinct points in
$[a,b] \times \{t_0\}$ touched by backward paths which also touch
$\R \times \{t_0 + t\}$. 
\erm

\brm
\label{rmk:ptchar}
As in~\cite{kn:TW}, space-time points $(x,t)$ can be characterized by the
number of locally disjoint paths $m_{\mbox{\rm\scriptsize in}}$ 
(resp., $m_{\mbox{\rm\scriptsize out}}$)
of the BW entering (resp., leaving) that point from earlier (resp., to later)
times. The corresponding dual BW characterization has 
$m^{\mbox{\rm\scriptsize dual}}_{\mbox{\rm\scriptsize in}}=
m_{\mbox{\rm\scriptsize out}}-1$ and
$m^{\mbox{\rm\scriptsize dual}}_{\mbox{\rm\scriptsize out}}=
m_{\mbox{\rm\scriptsize in}}+1$. Generic (e.g., deterministic) points have
$(m_{\mbox{\rm\scriptsize in}},m_{\mbox{\rm\scriptsize out}})=(0,1)$.
Almost surely, there are nongeneric points of type 
$(0,2),(0,3),(1,1),(1,2)$ and $(2,1)$ but no others. 
We note that as in~\cite{kn:TW}, ruling out points of higher type
uses improvements of (\ref{etabound}) for $k > 2$.
Type $(2,1)$ (resp., $(0,3)$) points are those where coalescing
(resp., dual coalescing) occurs. Type $(1,2)$ points are particularly
interesting in that the single incident path continues along exactly one 
of the two outward paths --- with the choice determined intrinsically
rather than by some convention.
\erm


\section{Convergence to the Brownian Web}
\setcounter{equation}{0}
\label{sec:conv}

Let $X_\d$ be an $(\h,\f_\h)$-valued random variable
indexed by $\d \in (0,1]$, with distribution $\mu_\d$. We present 
criteria sufficient to insure convergence in distribution 
as $\d\to0$ of $X_\d$ to the Brownian web $\bar W$,
in the setting where the $X_\d$'s have coalescing paths;
for simplicity, we will not present here more general criteria that
do not require the coalescing property. We next introduce the various
conditions on $\mu_\d$ which together will imply convergence.

\begin{figure}[!ht]
\begin{center}
\includegraphics[width=10cm]{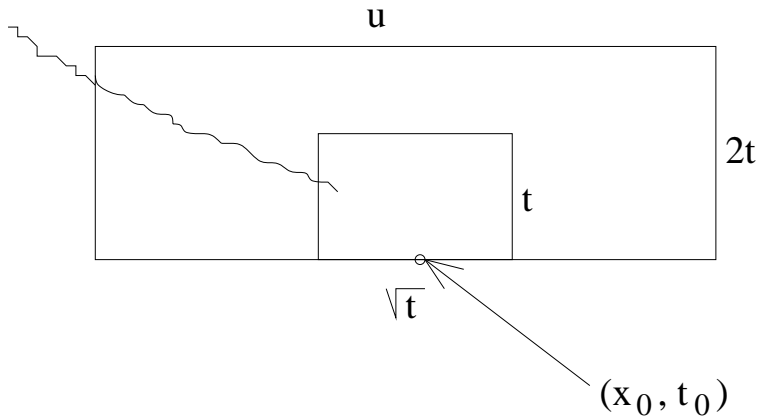}
\caption{Schematic diagram of a path causing the unlikely event
$A_{t,u}(x_0,t_0)$ to occur.} \label{figtight}
\label{int} 
\end{center}
\end{figure}

The first condition will guarantee tightness of the $\mu_\d$'s.
Let $R(x_0,t_0;u,t)$ denote the rectangle 
$[x_0-u/2,x_0+u/2]\times[t_0,t_0+t]$ in $\R^2$. We call 
$\{x_0\pm u/2\}\times[t_0,t_0+t]$ its right and left boundaries.
For $t>0,u>\sqrt t$, define $A_{t,u}(x_0,t_0)$ to be the event
(in $\f_\h$) that $K$ (in $\h$) contains a path touching both
$R(x_0,t_0;{\sqrt t},t)$ and (at a later time) the left or right 
boundary of the bigger rectangle $R(x_0,t_0;u,2t)$; see Figure~\ref{figtight}.
Our tightness condition is
\begin{equation}
\label{T1}
(T_1) \quad \tg(t,u)\equiv 
t^{-3/2}\sup_{\d>0}\,\,\sup_{x_0,t_0}\mu_\d(A_{t,u}(x_0,t_0))\to0\mbox{ as }t\to0\!+
\mbox{ for fixed } u>0\,.
\end{equation}

Our second condition will guarantee a weakened version of (i) in
Theorem~\ref{thm:char} (see Remark~\ref{rmk:cons}) for any limit
$\mu$ of $\mu_\d$. Let $\cd$ be any deterministic countable dense 
set of points in $\R^2$. The condition concerns the existence for
each $\d>0$ and $y\in\cd$ of measurable (on the probability space
of $X_\d$) single-path valued random variables 
$\theta^y_\d(\om)\in X_\d(\om):$

\vspace{.5cm}

$(I_1) $\quad There exist such $\theta^y_\d\in X_\d$ satisfying: 
for any deterministic $y_1,\ldots,y_m\in\cd$, 
${\theta^{y_1}_\d,\ldots,\theta^{y_m}_\d}$ converge in distribution
as $\d\to0$ to coalescing Brownian motions (with unit diffusion constant)
starting at $y_1,\ldots,y_m$.

\vspace{.5cm}

Our next two conditions will together guarantee (when $X_\d$ is
coalescing) a version of (ii) in Theorem~\ref{thm:char}
(see Remark~\ref{rmk:cons}).
For $-\infty<t_0<\infty$ and $0<t<\infty$, 

\begin{eqnarray}
\label{B1}
(B_1) &&\!\!\!\limsup_{\d \to 0}\sup_{a\in\R}
\mu_\d(\eta(t_0,t;a,a+\e)\geq2)\to0 \hbox{ as } \e\to0+;\\
\label{B2}
(B_2) &\e^{-1}&\!\!\!\limsup_{\d \to 0}\sup_{a\in\R}
\mu_\d(\eta(t_0,t;a,a+\e)\geq3)\to0 \hbox{ as } \e\to0+.
\end{eqnarray}

\bteo
\label{thm:conv}
Suppose $X_\d$ for $0<\d\leq1$ are $(\h,\f_\h)$-valued random variables 
with coalescing paths. If $T_1, I_1, B_1$ and $B_2$ all hold, then the 
distributions $\mu_\d$ of $X_\d$ converge weakly as $\d\to0$
to the distribution
$\mu_{\bar W}$ of the standard Brownian web.
\eteo

\noindent{\bf Sketch of Proof of Theorem~\ref{thm:conv}.}
We first explain why $T_1$ implies tightness. Let $g_{\d}(t,u)$
denote the sup over $x_0,t_0$ of $\mu_\d(A_{t,u})$
as in (\ref{T1}).
This represents an upper bound on the $\mu_\d$-probability
that there is 
some path $(f,t^*)$ passing through some point $(x',t')=(f(t'),t')$ in
the deterministic $\sqrt{t} \times t$ rectangle
$R(x_0,t_0;\sqrt{t},t)$ located at any $(x_0,t_0)$, 
such that for some $t''\in [t',t'+t]$
the spatial increment $|f(t'')-f(t')| \geq u$ even though
the time increment is $\le t$. Now taking a large $L \times T$
space-time rectangle centered at the origin and covering it with 
$O(LT/ t^{3/2})$
$\sqrt{t} \times t$ small rectangles, 
we see that $LT\,\tg(t,u)$ represents an
upper bound on the $\mu_\d$-probability (for any $\d$) that some path
has $|f(t'')-f(t')| \geq u$ while $t''-t' \le t$ with
$(f(t'),t')$ {\em anywhere} in the large rectangle. 
We next choose  sequences $u_n \to 0$, $L_n \to \infty$,
$T_n \to \infty$ and then $t_n \to 0$ sufficiently
rapidly that $L_nT_n\tg(t_n,u_n)$ is summable.
Now moving to the compactified space-time $\br^2$ (and
using the notation of (\ref{compactify})), it follows that there
are sequences $\phi_n,\psi_n \to 0$ so that for large
enough $n$, with $\mu_\d$-probability  close to one (for any $\delta$),
$|\Psi(t'')-\Psi(t')| \le \psi_n$ implies 
$\Phi(f(t''),t'')-\Phi(f(t'),t')|\, \le \,\phi_n$. This 
equicontinuity
with probability close to one (for any $\delta$) combined with
a version of the Arzel\`{a}-Ascoli theorem leads 
to the paths, as elements of $(\o,d)$,
belonging to a compact subset $K_{\tilde{\varepsilon}}$ of $(\o,d)$ with
$\mu_\d$-probability $\geq 1-\tilde{\varepsilon}$ 
(for any $\delta$), which implies tightness because the collection of
compact subsets of $K_{\tilde{\varepsilon}}$ is 
itself a compact set in $(\h,d_\h)$. 

Tightness implies that every subsequence of $\mu_\d$ has a
sub-subsequence 
converging weakly to some $\mu$. To complete the proof, we need
to show that any such $\mu$ equals $\mu_{\bar W}$. To do this,
we will show that $\mu$ satisifies the two characterization properties
of Theorem~\ref{thm:char}, as modified in Remark~\ref{rmk:cons}. 
Combining condition $I_1$ with convergence 
in distribution (along a subsequence)
of $X_\d$ to some $X$ distributed by $\mu$, we see that
we can realize $X$ on some probability space so that it contains
paths starting from the points of $\cd$ distributed as coalescing
Brownian motions. This is just the desired (weakened version of)
property(i) in
Theorem~\ref{thm:char}. Indeed, this shows that $X$ contains an
$X'$ that has the Brownian web distribution. 
 
To complete the proof we use 
conditions $B_1$ and $B_2$.
Note first that by limit/approximation arguments these two conditions
(without the $\limsup$ over $\d$) are valid with $\mu$ replacing
$\mu_\d$. For fixed $t_0,t,a,b$, we now consider $M+1$ equally spaced
points, $z_j = (a+j(b-a)/M,t_0)$ for $j=0,\dots,M$. For the random $X$
we will denote the counting variable $\eta(t_0,t;a,b)$ by $\eta$,
and the corresponding variable for $X'$ by $\eta'$. We also want
to count the number of points
on $\R\times\{t_0+t\}$ that are touched by paths
of $X$ that also touch $\{z_0,\dots,z_M\}$ and we will denote these
variables for $X$ and $X'$ by $\eta_M$ and $\eta'_M$. 
Of course, $\eta \geq \eta_M$ and $\eta \geq \eta'_M$.
By condition $B_1$ (for $\mu$) applied to small intervals about each
of the $z_j$'s, it follows that $\eta_M = \eta'_M$ almost surely.
Applying condition
$B_2$ (for $\mu$) to the $M$ spatial intervals $[z_{j-1},z_j]$
of length $\epsilon = (b-a)/M$,
and using the coalescing (or at least non-crossing)
property of $X$ that it inherits from
the $X_\d$'s, it follows that 
\begin{equation}
P(\eta > \eta'_M)\,=\, P(\eta > \eta_M)\to 0 \hbox{ as } M\to \infty.
\end{equation}
Thus $P(\eta > \eta')=0$
so that the distribution
of $\eta'$ is stochastically dominated by (and hence equal to)
the distribution of $\eta$. This gives the desired (modified version of)
property(ii) in Theorem~\ref{thm:char} and completes the proof.

\brm
\label{rmk:etaconv}
The arguments used in the proof of Theorem~\ref{thm:conv}
also show that under the same four conditions $T_1, I_1, B_1$ and $B_2$,
the counting random variables $\eta_\d(t_0,t;a,b)$ for $X_\d$ converge in
distribution to the Brownian web counting variable 
$\eta_{\bar W}(t_0,t;a,b)$.
If one also has uniform integrability as $\d \to 0$, then the means
converge to the mean of $\eta_{\bar W}$. 
\erm

To apply Theorem~\ref{thm:conv} to random walks, we
begin by precisely defining $Y$ (resp., $\tilde Y$),
the set of all discrete (resp., continuous) time coalescing random walks
on $\Z$. The sets of rescaled walks, $Y^{(\d)}$ and $\tilde Y^{(\d)}$,
are then obtained by
the usual rescaling of space by $\d$ and time by $\d^2$. 
The (main) paths of $Y$ are the discrete-time random walks
$\yo$, as described in the Introduction and 
shown in Figure~1, with $(y_0,s_0)=
(i_0,j_0) \in\Z \times \Z$ arbitrary except that $i_0+j_0$ must be even.
Each random walk path goes from $( i, j)$ to $(i\pm1,j+1)$
linearly. In addition to these, we add some boundary paths so that $Y$
will be a compact subset of $\o$. These are all the paths of the form
$(f,s_0)$ with $s_0\in \Z \cup \{-\infty,\infty\}$ and $f\equiv\infty$ or 
$f\equiv-\infty$. Note that for $s_0=-\infty$ there are two different
paths starting from the single point at $s_0=-\infty$ in $\br^2$.

The continuous time $\tilde Y$ can be defined similarly,
except that here $y_0$ is any $i_0\in\Z$ and $s_0$ is arbitrary
in $\R$. Continuous time walks are normally seen as jumping from 
$ i$ to $i\pm1$
at the times $T^{(i)}_k\in(-\infty,\infty)$ of a rate one Poisson
process. If the jump is, say, to $i+1$, then our polygonal path will
have a linear segment between $(i,T^{(i)}_k)$ and 
$(i+1,T^{(i+1)}_{k'})$, where $T^{(i+1)}_{k'}$ is the first Poisson
event at $i+1$ after $T^{(i)}_k$.
Furthermore, if $T^{(i_0)}_{k}<s_0<T^{(i_0)}_{k+1}$, then there will be
a constant segment in the path before the first nonconstant linear
segment. If $s_0=T^{(i_0)}_{k}$, then we take two paths: one with an
initial constant segment and one without.

\bteo
\label{thm:crwtobw}
Each of the collections of rescaled coalescing random
walk paths, $Y^{(\d)}$ (in discrete time) and 
$\tilde Y^{(\d)}$ (in continuous time) converges in distribution 
to the standard Brownian web as $\d \to 0$.
\eteo

\noindent{\bf Sketch of Proof of Theorem~\ref{thm:crwtobw}.}
By Theorem~\ref{thm:conv}, it suffices to verify conditions
$T_1, I_1, B_1$ and $B_2$. We will save the tightness condition
$T_1$ for last as it is the messiest to verify, at least in
the continuous time case of $\tilde Y^{(\d)}$. 

Condition $I_1$
is basically a consequence of the Donsker
invariance principle, as already noted in the Introduction.
Conditions $B_1$ and $B_2$ follow from the coalescing walks
version of the inequality of (\ref{etabound}), which is 
\begin{equation}
\label{walkbound}
\mu_\d(\eta(t_0,t;a,a+\e)\geq k)\, \leq \, 
[\mu_\d(\eta(t_0,t;a,a+\e)\geq2)]^{k-1}.
\end{equation}
Taking the sup over $a$ and the $\limsup$ over $\d$ and
using standard random walk arguments produces 
an upper bound of the from $C_k (\e/ \sqrt{t})^{k-1}$
which yields $B_1$ and $B_2$ as desired.

It remains to verify $T_1$. We will sketch the arguments for
the continuous time $\tilde Y^{(\d)}$; the discrete
time $Y^{(\d)}$ is easier and corresponds to a
portion of the continuous time arguments. As in the proof of
Theorem~\ref{thm:conv}, we denote by $g_{\d}(t,u)$
the sup over $x_0,t_0$ of $\mu_\d(A_{t,u})$,
where $\mu_\d$ now denotes the distribution of 
$\tilde Y^{(\d)}$. For the continuous
time case (and for $u \le 1$ and $\sqrt{t}$ much smaller than
$u$), we will obtain a $\delta$-independent
bound on $g_{\d}(t,u)$ that will yield $T_1$
by first obtaining, as we explain below, separate upper bounds in
three regions of $\d$-values that depend on $t,u$:
\begin{eqnarray}
\label{gbound1}
& C_1 \exp(-C_2 u/\sqrt{t}) & 
\hbox{for  }  D_0/\sqrt{t} \le \d^{-1}, \\
\label{gbound2}
& (2t/\d^2)^{u/(3\d)} 
& \hbox{for  }  6/u \le \d^{-1} \le D_0/\sqrt{t}, \\
\label{gbound3}
& C_3\,(t/u^2)^2 & \hbox{for  }  \d^{-1} \le 6/u.
\end{eqnarray}
Together, these bounds (with $D_0$ chosen appropriately) yield
\begin{equation}
\label{gbound4}
g_{\d}(t,u) \le C_4  t^2/u^4,
\end{equation}
which gives $T_1$ as desired. 

The first region of $\d$-values
corresponds to a spatial interval of
width $\sqrt{t}$ being multiple lattice spacings $\d$ wide and
a spatial interval of width $u$ being multiple $\sqrt{t}$-intervals
wide. The bound (\ref{gbound1}) comes about because the event 
$A_{t,u}$ is prevented if between the small rectangle and
both the left and right 
boundaries of the larger rectangle (see Figure~\ref{figtight}),
there is a random walk path that stays within some spatial 
$\sqrt{t}$-interval between times $t_0$ and $t_0 + 2t$. The second
region corresponds to two (or more) spatial lattice sites
between the small rectangle and the left (or right) 
boundary of the larger rectangle. The bound here comes from
preventing $A_{t,u}$ by
having a random walk path stay between two adjacent
spatial lattice sites between times $t_0$ and $t_0 + 2t$. The third bound
comes from preventing $A_{t,u}$ by not having the Poisson process occurences 
at adjacent spatial lattice sites,
$T^{(i)}_k$ and $T^{(i\pm 1)}_{k'}$, too close together in time.
This completes our sketch of the proof. 

\bigskip

\noindent  {\bf Acknowledgments.} Research partially supported by
FAPESP and CNPq (Brasil), MURST (Italia), and NSF (U.~S.~A.~). 
The authors thank S.~R.~S.~Varadhan for useful discussions. 


\end{document}